\documentclass{amsart}
\usepackage{times}

\begin{document}
\numberwithin{equation}{section}

\def\1#1{\overline{#1}}
\def\2#1{\widetilde{#1}}
\def\3#1{\widehat{#1}}
\def\4#1{\mathbb{#1}}
\def\5#1{\frak{#1}}
\def\6#1{{\mathcal{#1}}}

\newcommand{\de}{\partial}
\newcommand{\al}{\alpha}
\newcommand{\tr}{\widetilde{\rho}}
\newcommand{\tz}{\widetilde{\zeta}}
\newcommand{\tv}{\widetilde{\varphi}}
\newcommand{\hv}{\hat{\varphi}}
\newcommand{\tu}{\tilde{u}}
\newcommand{\tF}{\tilde{F}}
\newcommand{\debar}{\overline{\de}}
\newcommand{\Z}{\mathbb Z}
\newcommand{\C}{\mathbb C}
\newcommand{\Po}{\mathbb P}
\newcommand{\zbar}{\overline{z}}
\newcommand{\G}{\mathcal{G}}
\newcommand{\So}{\mathcal{S}}
\newcommand{\Ko}{\mathcal{K}}
\newcommand{\U}{\mathcal{U}}
\newcommand{\B}{\mathbb B}
\newcommand{\oB}{\overline{\mathbb B}}
\newcommand{\Cur}{\mathcal D}
\newcommand{\Dis}{\mathcal Dis}
\newcommand{\Levi}{\mathcal L}
\newcommand{\SP}{\mathcal SP}
\newcommand{\Sp}{\mathcal Q}
\newcommand{\A}{\mathcal O^{k+\alpha}(\overline{\mathbb D},\C^n)}
\newcommand{\CA}{\mathcal C^{k+\alpha}(\de{\mathbb D},\C^n)}
\newcommand{\Ma}{\mathcal M}
\newcommand{\Ac}{\mathcal O^{k+\alpha}(\overline{\mathbb D},\C^{n}\times\C^{n-1})}
\newcommand{\Acc}{\mathcal O^{k-1+\alpha}(\overline{\mathbb D},\C)}
\newcommand{\Acr}{\mathcal O^{k+\alpha}(\overline{\mathbb D},\R^{n})}
\newcommand{\Co}{\mathcal C}
\newcommand{\Hol}{{\sf Hol}(\mathbb H, \mathbb C)}
\newcommand{\Aut}{{\sf Aut}(\mathbb D)}
\newcommand{\D}{\mathbb D}
\newcommand{\oD}{\overline{\mathbb D}}
\newcommand{\oX}{\overline{X}}
\newcommand{\loc}{L^1_{\rm{loc}}}
\newcommand{\la}{\langle}
\newcommand{\ra}{\rangle}
\newcommand{\N}{\mathbb N}
\newcommand{\kd}{\kappa_D}
\newcommand{\Hr}{\mathbb H}
\newcommand{\ps}{{\sf Psh}}
\newcommand{\Hess}{{\sf Hess}}
\newcommand{\subh}{{\sf subh}}
\newcommand{\harm}{{\sf harm}}
\newcommand{\ph}{{\sf Ph}}
\newcommand{\tl}{\tilde{\lambda}}
\newcommand{\gdot}{\stackrel{\cdot}{g}}
\newcommand{\gddot}{\stackrel{\cdot\cdot}{g}}
\newcommand{\fdot}{\stackrel{\cdot}{f}}
\newcommand{\fddot}{\stackrel{\cdot\cdot}{f}}
\def\v{\varphi}
\def\Re{{\sf Re}\,}
\def\Im{{\sf Im}\,}

\title {A note on random holomorphic iteration in  convex domains}
\author[F. Bracci]{Filippo Bracci}
\address{Dipartimento di Matematica, Universit\`a di Roma  ``Tor Vergata'', Via della Ricerca Scientifica 1, 00133 Roma, Italy.}
\email{fbracci@mat.uniroma2.it} 

\begin{abstract} We introduce a geometric condition of Bloch type  which guarantees
that a subset of a bounded convex domain in several complex
variables is degenerate with respect to every iterated function
system. Furthermore we discuss the relations of such a Bloch type
condition with the analogous hyperbolic Lipschitz condition.
\end{abstract}

\maketitle

\def\Label#1{\label{#1}}


\def\cn{{\C^n}}
\def\cnn{{\C^{n'}}}
\def\ocn{\2{\C^n}}
\def\ocnn{\2{\C^{n'}}}
\def\je{{\6J}}
\def\jep{{\6J}_{p,p'}}


\def\dist{{\rm dist}}
\def\const{{\rm const}}
\def\rk{{\rm rank\,}}
\def\id{{\sf id}}
\def\aut{{\sf aut}}
\def\Aut{{\sf Aut}}
\def\CR{{\rm CR}}
\def\GL{{\sf GL}}
\def\Re{{\sf Re}\,}
\def\Im{{\sf Im}\,}
\def\U{{\sf U}}

\def\la{\langle}
\def\ra{\rangle}

\emergencystretch15pt \frenchspacing

\newtheorem{Thm}{Theorem}[section]
\newtheorem{Cor}[Thm]{Corollary}
\newtheorem{Pro}[Thm]{Proposition}
\newtheorem{Lem}[Thm]{Lemma}

\theoremstyle{definition}\newtheorem{Def}[Thm]{Definition}

\theoremstyle{remark}
\newtheorem{Rem}[Thm]{Remark}
\newtheorem{Exa}[Thm]{Example}
\newtheorem{Exs}[Thm]{Examples}

\def\bl{\begin{Lem}}
\def\el{\end{Lem}}
\def\bp{\begin{Pro}}
\def\ep{\end{Pro}}
\def\bt{\begin{Thm}}
\def\et{\end{Thm}}
\def\bc{\begin{Cor}}
\def\ec{\end{Cor}}
\def\bd{\begin{Def}}
\def\ed{\end{Def}}
\def\br{\begin{Rem}}
\def\er{\end{Rem}}
\def\be{\begin{Exa}}
\def\ee{\end{Exa}}
\def\bpf{\begin{proof}}
\def\epf{\end{proof}}
\def\ben{\begin{enumerate}}
\def\een{\end{enumerate}}

\section{Introduction}

Let $D\subset \C^n$ be a  domain. Let $\{f_j\}$ be a sequence of
holomorphic self-maps of $D$. Let $F_j:=f_j\circ\ldots \circ f_1$.
The sequence $\{F_j\}$ is called the {\sl holomorphic iterated
function system} (associated to $\{f_j\}$). Such systems are
encountered naturally in dynamical systems, continued fraction
theory and other areas of complex analysis. Given a holomorphic
iterated function system, one is interested in knowing its
asymptotical behavior, namely, to know the possible limits (in the
compact-open topology for instance) of the sequence. In general such
a question is rather difficult and one contents to know which
conditions guarantee that every limit of $\{F_j\}$ is constant (we
refer the reader to the papers \cite{BCMN} and \cite{Ke-La} and
bibliography therein).

More in detail, let $X\subset D$ be a subset  of $D$. We say that
the set $X$ is {\sl degenerate in $D$} if all the limits of any
holomorphic iterated function system  $\{F_j\}$ for which $f_j:D\to
X$ are constant.

In case $D=\D:=\{\zeta\in\D: |\zeta|<1\}$, degenerate subdomains are
completely characterized in terms of hyperbolic distance by Beardon,
Carne, Minda and Ng \cite{BCMN} and Keen and Lakic \cite{Ke-La}. To
state their results, we first introduce some terminology, as needed
for our later aims.

We denote by $k_D$ the Kobayashi distance of $D$ and by $\kappa_D$
the corresponding Kobayashi infinitesimal metric (for definition and
properties we refer to \cite{Ko}). Notice that for $D=\D$ then
$k_\D$ is nothing but the usual Poincar\'e distance.

Let $X\subset D$.  Let us denote by $R(X)$ its {\sl Bloch radius},
namely
\[
R(X)=\sup\{r\geq 0: B_D(z,r)\subset X\},
\]
where $B_D(w,r)=\{w\in D: k_D(w,z)<r\}$ is a Kobayashi ball of
center $z\in D$ and radius $r>0$.  A {\sl subset} $X\subset D$ is a
{\sl Bloch subset of $D$} if $R(X)<+\infty$.

Let now $Z\subseteq D$ be a {\sl subdomain} of $D$.  Then
$\kappa_Z(z;v)\geq \kappa_D(z;v)$ for all $z\in Z$ and $v\in T_zZ$.
The  {\sl hyperbolic Lipschitz constant} of $Z$ is defined as
\[
\mu_D(Z):=\sup_{z\in Z}\{\frac{\kappa_D(z;v)}{\kappa_Z(z;v)}:
(z,v)\in TZ, v\neq 0\}.
\]
The subdomain $Z$ is called a {\sl Lipschitz subdomain} of $D$ is
$\mu_D(Z)<1$.

In \cite{BCMN} Beardon, Carne, Minda and Ng proved the following
result:

\bt[Beardon-Carne-Minda-Ng]\label{Ng} Let $U\subset\D$ be a domain.
Then
\begin{enumerate}
\item $\tanh \frac{1}{2} R(U)\leq \mu_\D(U)\leq \tanh R(U)$. In
particular $U$ is a Bloch subdomain of $\D$ if and only if it is a
Lipschitz subdomain of $\D$.
\item If $U$ is a Lipschitz subdomain of $\D$ then $U$ is degenerate in $\D$.
\end{enumerate}\et

In \cite{Ke-La} Keen and Lakic showed that also the converse of (2)
holds:

\bt[Keen-Lakic]\label{KLac} Let $U\subset\D$ be a domain.  If $U$ is
degenerate in $\D$ then it is a Bloch subdomain of $\D$. \et

Thus the theorems of Beardon, Carne, Minda and Ng and of Keen and
Lakic completely characterize geometrically  degenerate subdomains
of $\D$.

In higher dimension the story is different. In \cite{BCMN} it is in
fact also proved:

\bt\label{Ng2} Let $D\subset\subset \C^n$ be a domain. If $X\subset
D$ is a Lipschitz subdomain of $D$ then $X$ is degenerate in $D$.
\et

However, as Beardon, Carne, Minda and Ng show, Bloch domains are not
necessarily Lipschitz domains in several dimensions, and  the
question of characterizing in a geometric flavor both Lipschitz and
degenerate subdomains of bounded domains in higher dimension is
open.

The aim of the present note is to present a   Bloch-type property
that guarantees a subset (not just an open subdomain) of a bounded
convex domain to be degenerate. To motivate our definitions and
results, we first look at the following example:

\be\label{esempio} Let $\{g_j\}$ be any sequence of holomorphic self
maps of $\D$ whose associated holomorphic iterated function system
 has some non-constant limit $g$ and whose image is
contained in some set $X'\subset\D$. Let us define $f_j:\B^n\to\B^n$
by $f_j(z_1,\ldots, z_n)=(g_j(z_1),0,\ldots,0)$. Then the
holomorphic iterated function system associated to $\{f_j\}$  has
image contained in $X= X'\times \{O\}$. The set $X$ is clearly a
Bloch subset of $\B^n$ because it contains no Kobayashi balls, but
the holomorphic iterated function system has a non-constant limit
$(g(z_1),0,\ldots, 0)$. Notice that however $X\cap (\D\times\{O\})$
cannot be a Bloch subdomain in $\D$ by Theorem \ref{KLac} and thus
it is not Lipschitz (as a subvariety) of $\B^n$ according to Theorem
\ref{Ng}. \ee

The previous example suggests that degenerate properties of a subset
should be related to Blochness properties of the intersection of
that subset with suitably chosen analytic discs. In order to make
this argument work (and to choose the right discs) we briefly recall
how a {\sl Lempert projection device} is defined (see \cite{Le},
\cite[Proposition 2.6.22]{Abate} and \cite[Theorem 4.8.12]{Ko} for
further details).

Let $D$ be a bounded   convex domain in $\C^n$  and let $z_0\in D$.
Given any point $z \in D$ there exists a
 {\sl complex geodesic} $\varphi: \D \to D$, {\em i.e.}, a
holomorphic isometry between $k_\D$  and $k_D$, such that
$\varphi(0)=z_0$  and $\varphi(t)=z$ for some $t\in (0,1)$. A
complex geodesic is also an infinitesimal isometry between the
Poincar\'e metric and the Kobayashi metric, and, given any point
$z\in D$ and nonzero direction $v\in T_zD$ there exists a  complex
geodesic containing $z$ and tangent to $v$ at $z$.

Moreover for any such a complex geodesic there exists a holomorphic
retraction $\rho_\v : D \to \varphi(\D)$ with affine fibers, {\em
i.e.} $\rho_\v$ is a holomorphic self-map of $D$ such that $\rho_\v
\circ \rho_\v =\rho_\v$, $\rho_\v(z)=z$ for any $z \in\varphi(\D)$
and $\rho_\v^{-1}(\v(\zeta))\cap D$ is the intersection of $D$ with
a complex hyperplane for all $\zeta\in\D$. We call such a $\rho_\v$
a {\em Lempert projection} associated to $\varphi$. We remark that
if $D$ is convex but not strongly convex then $\rho_\v$ is not
unique in general. For instance in the bidisc $\D\times \D$ the
complex geodesic $\D\ni\zeta\mapsto (\zeta, \zeta)$ has several
Lempert projections such as $\rho^1(z,w)=(z,z)$ and
$\rho^2(z,w)=(\frac{z+w}{2}, \frac{z+w}{2})$. However, if $D$ is
strongly convex then the Lempert projection (that is the one with
affine fibers) is unique (see \cite[Proposition 3.3]{BPT}).

Furthermore we let $\tr_\v:= \varphi^{-1} \circ \rho_\v:D\to \D$ and
call it the {\em left inverse } of $\varphi$, for $\tr_\v \circ
\varphi = {\sf id}_{\D}$. The triple $(\varphi, \rho_\v, \tr_\v)$ is
a so-called {\em Lempert projection device}.

\br For $D=\B^n$ the unit ball of $\C^n$ the image of the complex
geodesic through the points $z\neq w\in \overline{\B^n}$ is just the
{\sl one dimensional slice} $S_{z,w}:=\B^n\cap \{z+\zeta (z-w):
\zeta\in\C\}$. The Lempert projection is thus given by the
orthogonal projection of $\B^n$ onto $S_{z,w}$. \er

\bd Let $D\subset \C^n$ be a bounded   convex domain. We say that a
subset $X\subset D$ is {\sl 1-Bloch} in $D$  if there exists $C>0$
such that for any Lempert projection device $(\varphi, \rho_\v,
\tr_\v)$ the subset $\tr_\v(X)$ is contained in a Bloch subdomain
$U_\v$ of $\D$ with Bloch radius $R(U_\v)<C$. \ed

The main result of this note is:

\bt\label{main1} Let $D\subset \C^n$ be a bounded  convex domain.
Let $X\subset D$ be a subset of $D$ which is 1-Bloch in $D$. Then
$X$ is degenerate in $D$. \et

The proof of Theorem \ref{main1} is contained in section two. Such a
proof does not rely on any Lipschitzian property of 1-Bloch subsets
but, as we show in section three, a 1-Bloch {\sl subdomain} of $D$
is necessarily Lipschitz in $D$ (so that, in case $X$ is a
subdomain, Theorem \ref{main1} follows also from Theorem \ref{Ng2}).
In section three we also discuss of another natural Blochness
condition which is implied by the Lipschitz condition for
subdomains, giving some geometric hints on what a Lipschitz
subdomain in several complex variables looks like.

\section{Bloch, 1-Bloch and degenerate subsets}

In all the present section $D$ is   bounded   convex domain in
$\C^n$ and $X\subset D$ denotes a {\sl subset} of $D$. We begin with
the following simple observation:

\bp\label{chara} If $X$ is 1-Bloch than $X$ is a Bloch subset of
$D$.  \ep

\bpf Let $B_D(z_0,r)$ be any Kobayashi ball contained in $X$. Then
for any complex geodesic $\v:\D\to D$ such that $\v(0)=z_0$ it
follows that   $\v(B_\D(0,r))\subset B_D(z_0,r)$ and therefore
$\tr_\v(\v(\D)\cap X)$ contains the hyperbolic disc $B_\D(0,r)$.
Since $\tr_\v(\v(\D)\cap X)\subset \tr_\v(X)$ then $\tr_\v(X)$
contains the hyperbolic disc $B_\D(0,r)$.  Hence the Bloch radius of
any domain in $D$ containing $\tr_\v(X)$ must be greater than or
equal to $r$. Since by the very definition there exists a domain
$U_\v\subset \D$ containing $\tr_\v(X)$ with Bloch radius $\leq C$,
it follows that $r\leq C$ and hence $X$ has   Bloch radius $\leq C$,
and it is a Bloch subset of $D$. \epf

Notice that Example \ref{esempio} shows that the converse of
Proposition \ref{chara} is false in general.

In order to prove Theorem \ref{main1} we need a preliminary fact,
quite interesting by its own. First, we recall the following lemma
\cite[Lemma 3.1]{BCMN}

\bl\label{lemmaBMCN} If $g:\D\to U$ is holomorphic then
$k_\D(g(\zeta),g(\eta))\leq \mu(U)k_\D(\zeta,\eta)$ for all
$\zeta,\eta\in\D$. \el

Then we have

\bp \label{unioneBloch} Let $ W_j\subset \D$ with $j\in\N$ be Bloch
subdomains of $\D$. Assume that there exists $C>0$ such that the
Bloch radius $R(W_j)<C$ for all $j\in \N$. If $\{g_j\}$ is a
sequence of holomorphic self-maps of $\D$ such that
$g_j(\D)\subseteq W_j$ then the holomorphic iterated function system
$\{g_j\circ\ldots\circ g_1\}$ has only constant limits. \ep

\bpf By Theorem \ref{Ng} there exists $c<1$ such that the hyperbolic
Lipschitz constant $\mu(W_j)<c$ for all $j$. By Lemma
\ref{lemmaBMCN} it follows that $k_\D(g_j(\zeta), g_j(\eta))\leq c
k_\D(\zeta,\eta)$ for all $\zeta,\eta\in\D$ and $j\in\N$. Thus the
sequence $\{g_j\}$ is strictly decreasing with respect to $k_\D$
with uniform constant $c<1$ and our statement follows from a
Contraction Mapping Theorem  (see, \cite[Theorem 1.1]{BCMN}). \epf

Now we can prove Theorem \ref{main1}.

\begin{proof}[Proof of Theorem \ref{main1}]
Let $\{f_j\}$ be  a sequence of holomorphic self-maps of $D$ with
image contained in $X$. Let $F$ be a limit of the associated
holomorphic iterated system $\{F_j\}$ (where
$F_j:=f_j\circ\ldots\circ f_1$). Up to relabeling we can assume that
$\{F_j\}$ converges to $F$. Assume that $F$ is not constant and that
$F(z)\neq F(w)$ for some $z,w\in D$. Let $\v_1:\D\to D$ be a complex
geodesic such that $\v_1(0)=z$ and $\v_1(t_1)=w$ for some $0<t_1<1$.
Let $\v_j:\D\to D$ for $j\geq 2$ be a complex geodesic defined by
induction as follows: $\v_j(0)=f_{j-1}\circ \ldots\circ
f_1(z)=F_{j-1}(z)$ and $\v_j(t_j)=f_{j-1}\circ \ldots\circ
f_1(w)=F_{j-1}(w)$ for some $0<t_j<1$. Let $g_j:=\tr_{\v_{j+1}}\circ
f_j\circ \v_j$. Then $g_j:\D\to \D$ is holomorphic. By construction
$g_j\circ \ldots\circ g_1(0)=\tr_{\v_{j+1}}(F_j(z))$ while $g_j\circ
\ldots\circ g_1(t_1)=\tr_{\v_{j+1}}(F_j(w))$. Now, the family
$\{\v_{j}\}$ is a normal family in $D$. Let $\v:\D\to \overline{D}$
be one of its   limit. By continuity of the Kobayashi distance it
follows that either $\v$ is a complex geodesic or $\v(\D)\subset \de
D$. Since $\v_j(0)=F_{j-1}(z)\to F(z)$ as $j\to \infty$, it follows
that $\v$ is in fact a complex geodesic. Moreover, since
$\v_j(t_j)=F_j(w)\to F(w)$ as $j\to\infty$, it follows that there
exists $t\in (0,1)$ such that $t_j\to t$ as $j\to \infty$ and
$\v(t)=F(w)$. Up to re-labelling, we can assume that $\v$ is the
only limit of $\{\v_j\}$.

Next, we claim that the Lempert projections $\{\rho_{\v_{j}}\}$
converge (up to subsequences) in the compact-open topology of $D$ to
a Lempert projection $\rho_\v$ (showing that $\{\tr_{\v_j}\}$
converges to $\tr_\v$). To see this, let $\rho$ be any limit of the
normal family $\{\rho_{\v_{j}}\}$. Again, up to re-labeling, we can
assume that $\rho$ is the only limit. First we notice that, since
$\rho_{\v_j}(\v_j(\zeta))=\v_j(\zeta)$ for all $\zeta\in\D$ and
$j\in\N$, it follows that $\rho(\v(\zeta))=\v(\zeta)$ for all
$\zeta\in\D$. Thus, in particular, $\rho(D)\subset D$ and, since
$\rho_{\v_j}\circ \rho_{\v_j}=\rho_{\v_j}$ for any $j$, it follows
that $\rho\circ\rho=\rho$. Now we claim  that $\rho(D)=\v(\D)$. We
already know that $\v(\D)\subset\rho(D)$. Let $\rho(z)\in \rho(D)$.
Then there exists a sequence $\{z_{j}\}\subset D$ such that
$\rho_{j}(z_{j})\to \rho(z)$ as $j\to \infty$. But
$\rho_j(z_j)=\v_j(\zeta_j)$ for some $\zeta_j\in\D$ and $\zeta_j\to
\zeta\in \D$. Thus $\v_j(\zeta_j)\to\v(\zeta)=\rho(z)$ and hence
$\rho(z)\in\v(\D)$. Finally, it is    clear that the fibers of
$\rho$ are to be affine for those of every $\rho_{\v_j}$ are. Hence
$\rho$ is a Lempert projection associated to $\v$.

As a result, if $g$ is any limit of the holomorphic iterated
function system $\{g_j\}$ it follows that
$g(0)=\lim_{j\to\infty}\tr_{\v_{j+1}}(F_j(z))=\tr_\v(F(z))=0$ and
$g(t)=\lim_{j\to\infty}\tr_{\v_{j+1}}(F_j(w))=\tr_\v(F(w))=t$. Hence
the holomorphic iterated function system $\{g_j\}$ has a
non-constant limit. However, by the very definition of $1$-Bloch,
$g_j:\D\to W_j$ with $W_j=\tr_{\v_{j+1}}(f_j(\v_j(\D)))\subset
\tr_{\v_{j+1}}(X)\subset U_j$ where $U_j$ is a Bloch domain in $\D$
with Bloch radius bounded from above by some $C>0$ independently of
$j$. This contradicts Proposition~\ref{unioneBloch} and we are done.
\end{proof}

\section{c-Bloch, 1-Bloch and Lipschitz subdomains}

In all the present section $D$ is a bounded   convex domain in
$\C^n$ and $X\subset D$ denotes a {\sl subdomain} of $D$.

\bd  The  subdomain $X\subset D$ is {\sl c-Bloch} if there exists
$C>0$ such that for every Lempert projection device $(\varphi,
\rho_\v, \tr_\v)$  the (possibly empty) open set $\tr_\v(X\cap
\v(\D))$ is contained in a Bloch subdomain of $\D$ with Bloch radius
bounded from above by $C$. \ed

By the very definition we have

\bp A 1-Bloch subdomain $X\subset D$ is c-Bloch.\ep

The converse is however false in general as the following example
shows:

\be\label{esa2} Let $D=\B^2$. Let $E(e_1,R)$ denote a  {\sl
horosphere} (see, {\sl e.g.} \cite{Abate}) with center $e_1:=(1,0)$,
radius $R>0$ and pole $O$; namely
\[
E_{\B^2}(e_1,R)=\{z\in\B^2: |1-z_1|^2<R(1-\|z\|^2)\}.
\]
Let $X=E_{\B^2}(e_1,2)\setminus E_{\B^2}(e_1,1)$. The domain $X$ is
thus formed by the difference of two open complex ellipsoids tangent
each other to the point $(1,0)$. Its closure intersects the boundary
of $\B^2$ only at the point $(1,0)$. We claim that  $X$ is c-Bloch
but it is not 1-Bloch. Indeed the orthogonal projection of $X$ on
the complex geodesic $\v(\zeta)=(\zeta,0)$ is the horodisc
$E_\D(1,2)$ of $\D$ with center $1$ and radius $2$ which is not a
Bloch subdomain of $\D$, showing that $X$ is not 1-Bloch. To see
that $X$ is c-Bloch one can argue as follows. If $\v:\D\to D$ is any
complex geodesic whose closure contains $(1,0)$ then $\tr_\v(X\cap
\v(\D))$ is given by $E_\D(1,2)\setminus E_\D(1,1)$ (see
\cite[Proposition 2.7.8.(i)]{Abate}) which is a Bloch subdomain of
$\D$ (see \cite{BCMN}) with a fixed Bloch radius independent of
$\v$. As for the other complex geodesics, if $\v:\D\to D$ is any
complex geodesic whose closure does not contain $(1,0)$ then
$\tr_\v(X\cap \v(\D))$ is given by an annulus $A_\v$ in $\D$. Such
annuli $A_\v$ stay bounded from $\de\D$ independently of  $\v$
provided the $\v$'s stay away from a complex geodesic whose closure
contains $(1,0)$; while $A_\v$ ``degenerates'' into
$E_\D(1,2)\setminus E_\D(1,1)$ as $\v$ tends to a complex geodesic
whose closure contains $(1,0)$. Therefore the Bloch radius of
$\tr_\v(X\cap \v(\D))$ is bounded from above independently of $\v$.
\ee

We have the following relations among 1-Bloch, c-Bloch and Lipschitz
subdomain

\bp\label{main2} Let $D$ be a bounded  convex domain. Let $X\subset
D$ be a  subdomain.
\begin{enumerate}
\item If $X$  Lipschitz in $D$ then $X$ is c-Bloch in $D$.
\item If $X$ is 1-Bloch in $D$ then $X$ is Lipschitz in $D$.
\end{enumerate}
\ep

\begin{proof}
(1)   Let $\v:\D\to D$ be a complex geodesic such that $\v(\D)\cap
X\neq\emptyset$. The set $U_\v:=\tr_\v(X\cap \v(\D))$ is a (not
necessarily connected) domain in $\D$.  Since $\v|_{U_\v}: U_\v\to
X$ is holomorphic, by the monotonicity of the Kobayashi metric
$\kappa_{X}(\v(\zeta); d\v_\zeta (\xi))\leq
\kappa_{U_\v}(\zeta;\xi)$ for all $\zeta\in U_\v$ and $\xi\in
\C\setminus\{0\}=T_\zeta U_\v\setminus\{0\}$. Hence for all
$(\zeta,\xi)\in TU_\v$ we have
\[
  \frac{ \kappa_\D(\zeta; \xi)}{  \kappa_{U_\v}(\zeta;
 \xi)}=\frac{ \kappa_D(\v(\zeta); d\v_\zeta(
\xi))}{  \kappa_{U_\v}(\zeta;
 \xi)}\leq \frac{ \kappa_D(\v(\zeta); d\v_\zeta(
\xi))}{  \kappa_{X}(\v(\zeta); d\v_\zeta (\xi))}\leq c
\]
for some $c<1$ independent of $\v$ (because $X$ is Lipschitz in $D$
by hypothesis). Thus $U_\v$ is Lipschitz in $\D$ and by Theorem
\ref{Ng}.(1) it is a Bloch subdomain of $\D$ with Bloch radius
bounded from above by $2\tanh^{-1} c$. Hence $X$ is c-Bloch.

(2) Assume that $X$ is  1-Bloch. This means that there exists a
constant $C>0$ such that for all Lempert projection  devices
$(\varphi, \rho_\v, \tr_\v)$ it follows that the Bloch radius of
$\tr_\v(X)$ is less than or equal to $C$. In particular by
Theorem~\ref{Ng}.(1), the (possibly empty) open set $\tr_\v(X)$ is a
Lipschitz subdomain of $\D$ with hyperbolic Lipschitz constant
bounded from above by $c=\tanh C$.

Fix $z\in X$ and $v\in T_zD\setminus\{0\}$.  Let $\v:\D \to D$ be a
complex geodesic such that $\v(0)=z$ and $d\v_0(\xi)=v$ for some
$\xi\in \C$. Let $\tr_\v:D\to \D$ be the left inverse of $\v$. By
the monotonicity of the Kobayashi metric, considering the
holomorphic map $\tr_\v|_X:X\to \tr_\v(X)\subset \D$ and since
$d(\tr_\v)_{\v(\zeta)}\circ d\v_\zeta={\sf id}$, we have
\[
\kappa_{X}(z;v)\geq \kappa_{\tr_\v(X)}(\tr_\v(z);
d(\tr_\v)_z(v))=\kappa_{\tr_\v(X)}(0;\xi).
\]
Therefore, taking into account that
$\kappa_D(z;v)=\kappa_\D(0;\xi)$, it follows
\[
\frac{\kappa_D(z;v)}{\kappa_X(z;v)}\leq
\frac{\kappa_\D(0;\xi)}{\kappa_{\tr_\v(X)}(0;\xi)}\leq
\mu_\D(\tr_\v(X))\leq c,
\]
proving that $X$ is Lipschitz in $D$.\end{proof}

\br Proposition \ref{main2} gives a geometric necessary condition
(c-Blochness) for a subdomain to be Lipschitz, and then degenerate.
Such a condition is rather easy to be verified in simple domains
such as the unit ball $\B^n$ of $\C^n$.  We do not know whether such
a condition is also sufficient, namely, it is an open question if
c-Bloch implies Lipschitz. \er


\begin{thebibliography}{Co-MAc}
\bibitem{Abate} M. Abate, {\sl Iteration Theory of Holomorphic Maps on
Taut Manifolds}, Mediterranean Press, Rende, Cosenza, 1989.
\bibitem{BCMN} A. F. Beardon, T. K. Carne, D. Minda, T. W. Ng, {\sl
Random iteration of analytic maps}. Ergodic Th. Dyn. Systems, 24, 3,
(2004), 659-675.
\bibitem{Br} F. Bracci,  {\sl Dilatation and order of contact for holomorphic
self-maps of strongly convex domains}, Proc. London Math. Soc., 86,
1, (2003), 131-152.
\bibitem{BPT} F. Bracci, G. Patrizio, S. Trapani, {\sl The pluricomplex Poisson kernel for strongly convex
domains}, preprint 2005 (available on ArXiV).
\bibitem{Ke-La} L. Keen, N. Lakic, {\sl Random holomorphic
iterations and degenerate subdomains of the unit disk}. Proc. Amer.
Math. Soc. 134, 2, (2005) 371-378.
\bibitem{Ko} S. Kobayashi, {\sl Hyperbolic complex spaces}.
Springer, Grundlehren der mathematischen Wissenschaften 318.
\bibitem{Le} L. Lempert, {\sl La m\'etrique de Kobayashi et
la  representation des domaines sur la boule}. Bull. Soc. Math. Fr.
109  (1981), 427-474.
\end{thebibliography}
\end{document}